\documentclass[12pt]{article}
\usepackage{amssymb, amsmath, url, graphicx}

\def\3{\subset }
\def\4{\subseteq }
\def\<{\left<}
\def\>{\right>}

\def\bit{\begin{itemize}}
\def\eit{\end{itemize}}
\def\3{\subset }
\def\4{\subseteq }
\def\ov{\overline}

\def\0{\leqno}

\def\barr{\begin{array}}
\def\earr{\end{array}}
\def\dd{\displaystyle}

\def\Z{{\rlap{$\kern2pt{\rm Z}$}{\rm Z}\,}}
\def\bld#1#2{{\buildrel{#1}\over{#2}}}
\def\st#1#2{{\mathrel{\mathop{#2}\limits_{#1}}{}\!}}
\def\stb#1#2#3{{\st{{#1}}{\bld{{#2}}{#3}}{}\!}}
\def\xmare#1#2{\stb{#1}{#2}{\mbox{\Huge$\times$}}}

\def\frax{\dd\frac}

\title{\bf Cyclic factorization numbers\\ of finite groups}
\author{Marius T\u arn\u auceanu and Mihai-Silviu Lazorec}
\date{February 5, 2017}

\begin{document}

\maketitle

\begin{abstract}
In this  paper we introduce and study the concept of cyclic
factorization number of a finite group $G$. By using the
M\"{o}bius inversion formula and other methods involving the cyclic 
subgroup structure, this is explicitly computed for some
important classes of finite groups.
\end{abstract}

\noindent{\bf MSC (2010):} Primary 20D40; Secondary 20D60.

\noindent{\bf Key words:} cyclic factorization number, cyclic
subgroup commutativity degree, M\"{o}bius function.

\section{Introduction}

Let $G$ be a finite group, $L(G)$ be the subgroup lattice of $G$
and $L_1(G)$ be the poset of cyclic subgroups of $G$. One of the
most interesting positive integers that can be associated to $G$
is the \textit{fac\-to\-ri\-za\-tion number} $F_2(G)$. This counts
the number of all pairs $(H,K)\in L(G)^2$ satisfying $G=HK$ and
has been investigated in many recent papers, such as \cite{5,12}. On
the other hand, in our previous paper \cite{9} (see also
\cite{11}) we defined the \textit{subgroup commutativity degree}
$sd(G)$ of $G$ as the proportion of the number of ordered pairs
$(H,K)\in L(G)^2$ such that $HK=KH$ by $|L(G)|^2$. These two
quantities are closely connected. More precisely, we have
$$sd(G)=\dd\frac{1}{|L(G)|^2}\sum_{H\leq G} F_2(H)\0(1)$$and
$$F_2(G)=\dd\sum_{H\leq G} sd(H) |L(H)|^2\mu(H,G)\,,\0(2)$$by
applying the well-known M\"{o}bius inversion formula. We also
recall the following theorem due to P. Hall \cite{2} (see also
\cite{3}), that permits us to compute explicitly the M\"{o}bius
function of a finite $p$-group.

\bigskip\noindent{\bf Theorem 1.1.} {\it Let $G$ be a finite $p$-group of order $p^n$.
Then $\mu(G)=0$ unless $G$ is elementary abelian, in which case we
have $\mu(G)=(-1)^n p^{\,\binom{n}{2}}.$}
\bigskip

Then, in \cite{14}, we defined the \textit{cyclic subgroup
commutativity degree} $csd(G)$ of $G$ by replacing $L(G)$ with
$L_1(G)$ in the definition of $sd(G)$. This also has a
probabilistic significance, namely it measures the probability
that two cyclic subgroups of $G$ commute. So, it is naturally to
introduce a new positive integer that corresponds to $csd(G)$ in
the same way as $F_2(G)$ corresponds to $sd(G)$. In the following
we will denote by $CF_2(G)$ the number of all pairs $(H,K)\in
L_1(G)^2$ satisfying $G=HK$ and we will call it the \textit{cyclic
factorization number} of $G$. Its study is the purpose of the
current paper.
\bigskip

The paper is organized as follows. Some basic properties of cyclic
fac\-to\-ri\-za\-tion numbers are presented in Section 2. Section
3 deals with the computation of $CF_2(G)$ for some
classes of finite groups. In
the final section some further research directions and a list of
open problems are indicated.
\bigskip

Most of our notation is standard and will usually not be repeated
here. Elementary notions and results on groups can be found in
\cite{4,7}. For subgroup lattice concepts we refer the reader to
\cite{6,8,13}.

\section{Basic properties of cyclic factorization\\ numbers}

Let $G$ be a finite group. First of all, we remark that
$CF_2(G)\leq F_2(G)$, and we have equality if and only if $G$ is
cyclic. Consequently, for such a group $G$ the number $CF_2(G)$ is
given by \textbf{Theorem 1.1.} of \cite{5}.

\bigskip\noindent{\bf Proposition 2.1.} {\it Let $G$ be a finite
cyclic group of order $n=p_1^{\alpha_1}p_2^{\alpha_2}\cdots
p_k^{\alpha_k}$. Then
$$CF_2(G)=\Phi_1(n)=\prod_{i=1}^k (2\alpha_i+1)\,.$$}

Given two finite groups $G$ and $G'$, if $G\cong G'$ then
$CF_2(G)=CF_2(G')$. By \textbf{Proposition 2.1.} we infer that the
converse is not true. Also, we note that the weaker condition
$L(G)\cong L(G')$ does not imply $CF_2(G)=CF_2(G')$,
as shows the following elementary example.

\bigskip\noindent{\bf Example 2.2.} It is well-known that the subgroup
lattices of $G=\mathbb{Z}_3\times\mathbb{Z}_3$ and $G'=S_3$ are
isomorphic. On the other hand, we can easily check that
$CF_2(G)=12\neq 6=CF_2(G')$.
\bigskip

By a direct calculation, one obtains
$$CF_2(S_3\times\mathbb{Z}_2)=12\neq 18=CF_2(S_3)CF_2(\mathbb{Z}_2)$$and
therefore in general we don't have $$CF_2(G\times
G')=CF_2(G)CF_2(G').$$A sufficient condition in order to this
equality holds is that $G$ and $G'$ be of coprime orders. This
remark can naturally be extended to arbitrary finite direct
products.

\bigskip\noindent{\bf Proposition 2.3.} {\it Let
$(G_i)_{i=\overline{1,k}}$ be a family of finite groups having
coprime orders. Then
$$CF_2(\xmare{i=1}k G_i)=\prod_{i=1}^k CF_2(G_i)\,.$$}
\smallskip

The following immediate consequence of \textbf{Proposition 2.3.} shows that
the computation of $CF_2(G)$ for a finite nilpotent
group $G$ can be reduced to $p$-groups.

\bigskip\noindent{\bf Corollary 2.4.} {\it If $G$ is a finite
nilpotent group and $(G_i)_{i=\ov{1,k}}$ are the Sylow subgroups
of $G$, then $$CF_2(G)=\prod_{i=1}^k CF_2(G_i)\,.$$}

\bigskip\noindent{\bf Remarks 2.5.}
\begin{itemize}
\item[{\rm a)}]The conclusion of \textbf{Corollary 2.4.} fails if $G$ is not nilpotent.
For example, we have $CF_2(S_4)=CF_2(A_4)=0$ even if $CF_2(S)\neq
0$ for every Sylow subgroup $S$ of $S_4$ and $A_4$, respectively.
Notice also that more can be said about such groups, namely
$CF_2(S_n)=CF_2(A_n)=0$ for all $n\geq 4$.
\item[{\rm b)}] Let $p\geq 3$ be a prime and $G$ be a finite
$p$-group. It is well-known (see e.g. \cite{4}, I) that $G$ can
be written as a product of two cyclic subgroups if and only if it
is metacyclic. This can be reformulated in the following nice way:
"A finite $p$-group with $p\geq 3$ is metacyclic if and only if its
cyclic factorization number is non-zero".
\end{itemize}

Finally, we observe that the connections between $CF_2(G)$ and
$csd(G)$ are similar with (1) and (2), more exactly
$$csd(G)=\dd\frac{1}{|L_1(G)|^2}\sum_{H\leq G} CF_2(H)\0(3)$$and
$$CF_2(G)=\dd\sum_{H\leq G} csd(H) |L_1(H)|^2\mu(H,G)\,.\0(4)$$The
equality (4) is the main ingredient that will be used in Section 3
to calculate the cyclic factorization numbers of certain finite
groups.

\section{Cyclic factorization numbers for some classes of finite groups}

As we already have seen, the computation of cyclic factorization
numbers of finite abelian groups is reduced to $p$-groups. By the
fundamental theorem of finitely generated abelian groups, such a
group is of type $$G\cong \xmare{i=1}k
\mathbb{Z}_{p^{\alpha_i}}\,,$$where
$1\leq\alpha_1\leq\alpha_2\leq...\leq\alpha_k$. Recall that $G$
possesses a unique maximal elementary abelian subgroup
$M\cong\mathbb{Z}_p^k$. Moreover, all elementary abelian subgroups
of $G$ are contained in $M$. Notice also that we have $csd(H)=1$,
for every $H\leq G$. We are now able to compute explicitly
$CF_2(G)$.

\bigskip\noindent{\bf Theorem 3.1.} {\it The cyclic factorization number
of the finite abelian $p$-group $G\cong \xmare{i=1}k
\mathbb{Z}_{p^{\alpha_i}}$,
$1\leq\alpha_1\leq\alpha_2\leq...\leq\alpha_k$, is given by the
following equality:
$$CF_2(G)=\left\{\barr{lll}
2\alpha_1{+}1,& \mbox{if } k=1\\
\\
p^{2\alpha_1{-}1}\left[(2\alpha_2{-}2\alpha_1{+}1)p{-}2\alpha_2{+}2\alpha_1{+}1\right],& \mbox{if } k=2\\
\\
0,& \mbox{if } k\geq 3\,.\earr\right.$$}

\noindent{\bf Proof.} By using (4) and \textbf{Theorem 1.1.}, one obtains
$$\hspace{-15mm}CF_2(G)=\dd\sum_{H\leq G} |L_1(H)|^2\mu(H,G)=\dd\sum_{H\leq G} |L_1(G/H)|^2\mu(H)=\0(5)$$
$$\hspace{23mm}=\dd\sum_{H\leq M} |L_1(G/H)|^2\mu(H)=\dd\sum_{i=0}^k \,\sum_{H\leq M,\,|H|=p^i} \hspace{-2mm}|L_1(G/H)|^2 (-1)^i
p^{\,\binom{i}{2}}\,.$$For $k=1$ we have
$$CF_2(G)=|L_1(\mathbb{Z}_{p^{\alpha_1}})|^2-|L_1(\mathbb{Z}_{p^{\alpha_1-1}})|^2=(\alpha_1+1)^2-\alpha_1^2=2\alpha_1+1\,,$$while
for $k=2$ we have
$$CF_2(G)=|L_1(\mathbb{Z}_{p^{\alpha_1}}\times\mathbb{Z}_{p^{\alpha_2}})|^2-p|L_1(\mathbb{Z}_{p^{\alpha_1-1}}\times\mathbb{Z}_{p^{\alpha_2}})|^2-|L_1(\mathbb{Z}_{p^{\alpha_1}}\times\mathbb{Z}_{p^{\alpha_2-1}})|^2+$$
$$\hspace{15mm}+p|L_1(\mathbb{Z}_{p^{\alpha_1-1}}\times\mathbb{Z}_{p^{\alpha_2-1}})|^2=p^{2\alpha_1{-}1}\left[(2\alpha_2{-}2\alpha_1{+}1)p{-}2\alpha_2{+}2\alpha_1{+}1\right]$$by
Theorem 4.2 of \cite{10}.
\smallskip

In the case $k\geq 3$ the desired equality $CF_2(G)=0$ follows by
a direct calculation from (5) and \textbf{Theorem 4.3.} of \cite{10}, or by observing
that $G$ cannot be generated by two elements (and consequently cannot
be written as a product of two cyclic subgroups). This completes the proof.
\hfill\rule{1,5mm}{1,5mm}
\bigskip

In the second part of this section we will focus on computing the
cyclic factorization number of other classes of finite groups starting with 
the dihedral groups
$$D_{2n}=\langle x,y\mid x^n=y^2=1,\ yxy=x^{-1}\rangle,\, n\geq 3\,.$$The subgroup
structure of $D_{2n}$ is the following: for every divisor $d$ or
$n$, $D_{2n}$ possesses a subgroup isomorphic to $\mathbb{Z}_d$,
namely $H_d=\langle x^{\frac{n}{d}}\rangle$, and $\frac{n}{d}$
subgroups isomorphic to $D_{2d}$, namely $K_d^i=\langle
x^{\frac{n}{d}},x^{i-1}y\rangle$, $i=1,2,...,\frac{n}{d}\,.$ We
easily infer that $$|L_1(H_d)|=\tau(d) \mbox{ and }
|L_1(K_d^i)|=\tau(d)+d\,,\0(6)$$where $\tau(d)$ denotes the number
of divisors of $d$. On the other hand, we know that
$$csd(H_d){=}1 \mbox{ and } csd(K_d^i){=}\hspace{-1mm}\left\{\barr{lll}
\hspace{-1mm}\frax{\tau(d)(\tau(d){+}d){+}d(\tau(d){+}1)}{(\tau(d){+}d)^2},& \hspace{-1mm}d\equiv 1 \hspace{0,5mm}({\rm mod}\hspace{0,5mm} 2)\\
\\
\hspace{-1mm}\frax{\tau(d)(\tau(d){+}d){+}d(\tau(d){+}2)}{(\tau(d){+}d)^2},& \hspace{-1mm}d\equiv
0 \hspace{0,5mm}({\rm mod}\hspace{0,5mm} 2)\earr\right.\0(7)$$by
\textbf{Theorem 3.2.1.} of \cite{14}. Also, the values of the M\"{o}bius function
associated to $L(D_{2n})$ have been determined in \textbf{Lemma 19} of \cite{1}:
$$\mu(H_d,D_{2n})=\mu(D_{2\frac{n}{d}})=-\frac{n}{d}\mu(\frac{n}{d})\mbox{ and }\mu(K_d^i,D_{2n})=\mu(\mathbb{Z}_{\frac{n}{d}})=\mu(\frac{n}{d})\0(8)$$for
all divisors $d$ of $n$ and all $i=1,2,...,\frac{n}{d}\,.$ We are now able to complete the computation
of $CF_2(D_{2n})$.

\bigskip\noindent{\bf Theorem 3.2.} {\it The cyclic factorization number
of the dihedral group $D_{2n}$, $n\geq 3$, is given by the
following equality:
$$CF_2(D_{2n})=2n\,.$$}

\noindent{\bf Proof.} By using (6)-(8) in (4), one obtains
$$\hspace{-20mm}CF_2(D_{2n})=\dd\sum_{H\leq D_{2n}} csd(H) |L_1(H)|^2\mu(H,D_{2n})=$$
$$=\sum_{d|n} csd(H_d) |L_1(H_d)|^2\mu(H_d,D_{2n})+\sum_{d|n} \frac{n}{d}\, csd(K_d^1) |L_1(K_d^1)|^2\mu(K_d^1,D_{2n})=$$
$$=\sum_{d|n} \tau(d)^2\left(-\frac{n}{d}\mu(\frac{n}{d})\right)+\sum_{d|n} \frac{n}{d}\, \left[\tau(d)(\tau(d)+d)+d(\tau(d)+i_d)\right]\mu(\frac{n}{d})=$$
$$=2n\sum_{d|n} \tau(d)\mu(\frac{n}{d})+n\sum_{d|n} i_d\mu(\frac{n}{d}),$$where $$i_d=\left\{\barr{lll}
1,& d\equiv 1 \hspace{1mm}({\rm mod}\hspace{1mm} 2)\\
\\
2,& d\equiv 0 \hspace{1mm}({\rm mod}\hspace{1mm} 2)\earr\right.$$We observe that
$$\sum_{d|n} \tau(d)\mu(\frac{n}{d})=1 \mbox{ and } \sum_{d|n} i_d\mu(\frac{n}{d})=0$$and consequently
$$CF_2(D_{2n})=2n\,,$$as desired. \hfill\rule{1,5mm}{1,5mm}

\bigskip\noindent{\bf Remark 3.3.} An alternative and more facile way of computing $CF_2(D_{2n})$ uses its definition.
Clearly, we have $$L_1(D_{2n})=\{H_d \mid d|n\}\cup\{K_1^i \mid i=1,2,...,n\}.$$We
easily infer that $D_{2n}$ is the product of two cyclic subgroups if and only if one of them is $H_n$ and the other is of type
$K_1^i$ with $i=1,2,...,n$. Hence $CF_2(D_{2n})=2n$.\\

The computation of the factorization numbers of the dihedral group $D_{2n}$ and of the abelian $p$-group $G\cong \xmare{i=1}k
\mathbb{Z}_{p^{\alpha_i}}$,
$1\leq\alpha_1\leq\alpha_2\leq...\leq\alpha_k$, will be extremely helpful 
in our study. These results will be used to find the cyclic factorization numbers of the following 5 groups:\\
- the generalized quaternion group $$Q_{2^n}=\langle x,y \ | \ x^{2^{n-1}}=y^4=1, yxy^{-1}=x^{-1}\rangle, n\geq 3$$
- the quasi-dihedral group $$S_{2^n}=\langle x,y \ | \ x^{2^{n-1}}=y^2=1, yxy=x^{2^{n-2}-1}\rangle, n\geq 4$$
- the modular $p$-group $$M(p^n)=\langle x,y \ | \ x^{p^{n-1}}=y^p=1, y^{-1}xy=x^{p^{n-2}+1}\rangle, n\geq 3$$
- the dicyclic group $$Dic_{4n}=\langle a,\gamma \ | \ a^{2n}=1, \gamma^2=a^n, a^{\gamma}=a^{-1}\rangle,$$
- the generalized dicyclic group $$Dic_{4n}(A)=\langle A, \gamma \ | \ \gamma^4=e, \gamma^2\in A\setminus \lbrace e\rbrace, g^{\gamma}=g^{-1}, \forall g\in A\rangle,$$ where $A$ is an abelian group of order $2n$ isomorphic to $\mathbb{Z}_2\times \mathbb{Z}_n.$\\

As we saw in the alternative proof concerning the cyclic factorization number of the dihedral group $D_{2n}$, it is important to know the cyclic subgroup structure of the group we study. Hence, it is useful to recall that the poset of cyclic subgroups of $Q_{2^n}$ is 
$$L_1(Q_{2^n})=L(\langle x\rangle)\cup \lbrace \langle x^ky\rangle \ | \ k=\overline{0,2^{n-2}-1}\rbrace.$$
We are ready to compute the factorization number of the generalized quaternion group.\\

\noindent\textbf{Theorem 3.4.} \textit{The cyclic factorization number of the generalized quaternion group $Q_{2^n}$, $n\geq 3$, is given by
$$CF_2(Q_{2^n})=\begin{cases} 6, &\mbox{if } n=3 \\  2^{n-1}, &\mbox{if } n\geq 4 \end{cases}.$$}\\

\noindent\textbf{Proof.}  
By inspecting the subgroup lattice of $Q_8$, we infer that $CF_2(Q_8)=6$. Let $n\geq 4$ be a positive integer. The unique minimal subgroup of the generalized quaternion group $Q_{2^n}$ is its center $Z(Q_{2^n})=\langle x^{2^{n-2}}\rangle$. Also, it is known that 
$$\frac{Q_{2^n}}{Z(Q_{2^n})}\cong D_{2^{n-1}}.$$
Let $H$ and $K$ be two nontrivial cyclic subgroups of $Q_{2^n}$ such that $HK=Q_{2^n}$. Then the two subgroups contain $Z(Q_{2^n})$ and by the above isomorphism, the cyclic factorization $(H,K)$ of $Q_{2n}$ corresponds to the cyclic factorization $(\frac{H}{Z(Q_{2^n})},\frac{K}{Z(Q_{2^n})})$ of $D_{2^{n-1}}$. Also, it is important to remark that by quotiening non-cyclic subgroups of $Q_{2^n}$ by the center of the group, we do not obtain corresponding cyclic subgroups of $D_{2^{n-1}}$. In fact, by this procedure, we always get a dihedral group or a direct product of two abelian groups, namely $\mathbb{Z}_2 \times \mathbb{Z}_2$. With this being checked, we are sure that there are no cyclic factorizations of $D_{2^{n-1}}$ that could correspond to non-cyclic factorization of $Q_{2^n}$. We add that it is obvious that we can not build another cyclic factorization of the generalized quaternion group using its trivial subgroup as a factor. Hence $CF_2(Q_{2^n})=CF_2(D_{2^{n-1}})=2^{n-1}$ and our proof is complete.\hfill\rule{1,5mm}{1,5mm}\\

We will see that it may happen that by quotiening we can find a non-cyclic factorization of a group which corresponds to a cyclic factorization of its factor. The class of quasi-dihedral groups provides such an example. First of all, we recall that the cyclic subgroups of $S_{2^n}$ are illustrated by the poset
$$L_1(S_{2^n})=L(\langle x\rangle)\cup \lbrace \langle x^{2k}y\rangle \ | \ k=\overline{0,2^{n-2}-1}\rbrace\cup \lbrace \langle x^{2k+1}y\rangle \ | \ k=\overline{0,2^{n-3}-1} \rbrace.$$
Using this information, we can find the factorization number of the quasi-dihedral group.\\

\noindent\textbf{Theorem 3.5.} \textit{The cyclic factorization number of the quasi-dihedral group $S_{2^n}$, $n\geq 4$, is given by
$$CF_2(S_{2^n})=3\cdot 2^{n-2}.$$}\\

\noindent\textbf{Proof.} 
Besides the center $Z(S_{2^n})=\langle x^{2^{n-2}}\rangle$, the quasi-dihedral group has other  minimal subgroups contained in the set $M=\lbrace \langle x^{2k}y\rangle \ | \ k=\overline{0,2^{n-2}-1}\rbrace$. Moreover any nontrivial subgroup of $S_{2n}$ contains the center or it is contained in the set $M$. Let $H$ and $K$ be two cyclic subgroups such that $Z(S_{2^n})\subseteq H$, $Z(S_{2^n})\subseteq K$ and $HK=S_{2^n}$. In this case, to count the cyclic factorizations of $S_{2^n}$, we will use the isomorphism $$\frac{S_{2^n}}{Z(S_{2^n})}\cong D_{2^{n-1}}.$$
This result implies that a cyclic factorization $(H,K)$ of $S_{2^n}$ corresponds to a cyclic factorization of $D_{2^{n-1}}$. However, the quasi-dihedral group $S_{2^n}$ contains $2^{n-3}$ subgroups isomorphic to $\mathbb{Z}_2\times \mathbb{Z}_2$ and their images through the above isomorphism are cyclic subgroups of $D_{2^{n-1}}$ isomorphic to $\mathbb{Z}_2$. Each such subgroup is not contained in the maximal cyclic subgroup of $D_{2^{n-1}}$ and this leads us to $2\cdot 2^{n-3}=2^{n-2}$ cyclic factorizations of $D_{2^{n-1}}$ which correspond to non-cyclic factorizations of $S_{2^n}$. We remark that for other non-cyclic subgroups of the quasi-dihedral group $S_{2^n}$, we obtain non-cyclic subgroups of the dihedral group $D_{2^{n-1}}$ after quotiening by $Z(S_{2^n})$.
It is clear that we can not build a cyclic factorization of $S_{2^n}$ using two subgroups contained in $M$. Hence, the remaining option is to find cyclic factorizations $(H,K)$, where $H\in M$ and $K$ contains $Z(S_{2^n})$. Since, in this case, $H\cap K$ is trivial, we have only one choice for $K$, this being $K\cong \mathbb{Z}_{2^{n-1}}$, which is the only maximal cyclic subgroup of $S_{2n}$. By a simple computation, we get
$$CF_2(S_{2^n})=CF_2(D_{2^{n-1}})-2^{n-2}+2^{n-1}=3\cdot 2^{n-2},$$
as desired.\hfill\rule{1,5mm}{1,5mm}\\    

The class of modular $p$-groups $M(p^n)$ provided interesting probabilistic aspects such as $sd(M(p^n))=csd(M(p^n))=1$. The next purpose of this paper is to find an explicit formula for the cyclic factorization numbers of these groups.\\    

\noindent\textbf{Theorem 3.6.} \textit{The cyclic factorization number of the modular p-group $M_{p^n}$, $n\geq 3$, is given by
$$CF_2(M(p^n))=\begin{cases} 8, &\mbox{if } p=2,n=3 \\ p[(2n-4)(p-1)-p+3]+2p(p-1), &\mbox{if } p\geq 2,n\geq 4 \end{cases}.$$}\\

\noindent\textbf{Proof.} Using \textbf{Theorem 3.2.} and the isomorphism $M(8)\cong D_8$, we infer that $CF_2(M(8))=8$. Let $p\geq 2$ and $n\geq 4$ be a prime number and a positive integer, respectively. By inspecting the lattice of subgroups of $M(p^n)$, we remark that this group has $p+1$ minimal subgroups isomorphic to $\mathbb{Z}_p$, one of them being the commutator subgroup $D(M(p^n))=\langle x^{p^{n-2}}\rangle$ which is contained in each subgroup of higher order. Again, we will use an isomorphism to count the cyclic factorization number of the modular $p$-group. Indeed, it is known that 
$$\frac{M(p^n)}{D(M(p^n))}\cong \mathbb{Z}_p \times \mathbb{Z}_{p^{n-2}}.$$
Let $(H,K)$ be a cyclic factorization of $M(p^n)$ such that $D(M(p^n))\subseteq H$ and $D(M(p^n))\subseteq K$. Through the above isomorphism, there is a corresponding cyclic factorization of $\mathbb{Z}_p \times \mathbb{Z}_{p^{n-2}}$, namely $(\frac{H}{D(M(p^n))}, \frac{K}{D(M(p^n))})$. There are $n-2$ non-cyclic  proper subgroups of $M(p^n)$, these being isomorphic to $\mathbb{Z}_p\times \mathbb{Z}_{p^i}$, where $i=\overline{1,n-2}$. Quotiening these subgroups by $D(M(p^n))$, we obtain only one cyclic subgroup isomorphic to a minimal subgroup of $\mathbb{Z}_p\times \mathbb{Z}_{p^{n-2}}$. This subgroup can be used to obtain $2p$ cyclic factorizations of $\mathbb{Z}_p\times \mathbb{Z}_{p^{n-2}}$ since it is not contained in any maximal cyclic subgroup of this group. Hence, there are $2p$ cyclic factorizations of $\mathbb{Z}_p\times \mathbb{Z}_{p^{n-2}}$ corresponding to non-cyclic factorizations of $M(p^n)$. The only method left to find other cyclic factorizations $(H,K)$ of $M(p^n)$ is by using the $p$ minimal subgroups different from $D(M(p^n))$ as a factor and the subgroups that contain $D(M(p^n))$ as the other factor. In this case $H\cap K$ is trivial, so the second factor must be one of the $p$ maximal cyclic subgroups of $M(p^n)$ which are isomorphic to $\mathbb{Z}_{p^{n-1}}$, namely $\langle x\rangle$ and $\langle x^iy\rangle$, where $i=\overline{1,p-1}$. Adding up the numbers and using \textbf{Theorem 3.1.}, we have 
$$CF_2(M(p^n))=CF_2(\mathbb{Z}_p\times \mathbb{Z}_{p^{n-2}})-2p+2p^2=p[(2n-4)(p-1)-p+3]+2p(p-1),$$ as stated above. \hfill\rule{1,5mm}{1,5mm}\\

A simplified formula is obtained for the modular 2-groups.\\

\noindent\textbf{Corollary 3.7.} \textit{The cyclic factorization number of the modular 2-group $M(2^n)$, $n\geq 4$, is given by $$CF_2(M(2^n))=4n-2.$$}\\

Some probabilistic aspects of the dicyclic groups were studied in \cite{15}. We know that   
for each divisor $r$ of $2n$, the dicyclic group $Dic_{4n}$ has only one subgroup isomorphic to $\mathbb{Z}_r$, namely $\langle a^{\frac{2n}{r}}\rangle$. Also, for each divisor $s$ of $n$, $Dic_{4n}$ possesses $\frac{n}{s}$ subgroups isomorphic to $Dic_{4s}$, namely $\langle a^{\frac{n}{s}}, a^{i-1}\gamma\rangle$, where $i=\overline{1,\frac{n}{s}}$. The poset of cyclic subgroups is 
$$L_1(Dic_{4n})=\lbrace \langle a^{\frac{2n}{r}}\rangle \ | \ r|2n\rbrace \cup \lbrace \langle a^{i-1}\gamma\rangle \ | \ i=\overline{1,n}\rbrace.$$
We mention that, by definition, $Dic_4\cong \mathbb{Z}_4$. Also, it is easy to see that $Dic_8\cong Q_8$. Then $CF_2(Dic_4)=5$ and $CF_2(Dic_8)=6.$ Finding the factorization numbers for all other dicyclic groups is our next purpose. 
\\

\noindent\textbf{Theorem 3.8.} \textit{The cyclic factorization number of the dicyclic group $Dic_{4n}$, $n\geq 3$, is given by
$$CF_2(Dic_{4n})=\begin{cases} 4n, &\mbox{if } n\equiv 1(mod \ 2) \\ 2n, &\mbox{if } n\equiv 0(mod \ 2) \end{cases}.$$}\\

\noindent\textbf{Proof.}  It is known that $Z(Dic_{4n})=\langle a^n\rangle$ which is isomorphic to $\mathbb{Z}_2$ and the following isomorphism holds
$$\frac{Dic_{4n}}{Z(Dic_{4n})}\cong D_{2n}.$$
We remark that by quotiening non-cyclic subgroups of $Dic_{4n}$ we do not obtain cyclic subgroups of $D_{2n}$, so there will be no need to substract a number of factorizations from $CF_2(D_{2n})$ to find $CF_2(Dic_{4n})$, like we did in the case of other classes of groups. If $n=2^m$, where $m\geq 2$ is a positive integer, since $Dic_{4n}\cong Q_{4n}$, we have $CF_2(Dic_{4n})=CF_2(Q_{4n})=2n$ according to \textbf{Theorem 3.3.}. If $n$ is even and is not a power of 2, a nontrivial subgroup $H$ of $Dic_{4n}$ contains the center of the group or it is a subgroup of $\mathbb{Z}_{2n}$ of odd order. Let $(H,K)$ be a cyclic factorization of $Dic_{4n}$ such that $Z(Dic_{4n})\subseteq H$ and $Z(Dic_{4n})\subseteq K$. Then there is a corresponding cyclic factorization of the dihedral group $D_{2n}$. If $H$ is a nontrivial cyclic subgroup of odd order contained in the maximal subgroup $\mathbb{Z}_{2n}$ of $Dic_{4n}$, then we can use it to form a cyclic factorization if the other factor, $K$, is a cyclic subgroup which is not contained in $\mathbb{Z}_{2n}$. Looking at the poset $L_1(Dic_{4n})$, the only possible choices for $K$ are the cyclic subgroups of order 4 from the set $\lbrace \langle a^{i-1}\gamma\rangle \ | \ i=\overline{1,n}\rbrace$. This leads us to $|H|=n$, which is false since $n$ is even and $|H|$ is odd. Hence, if $n\equiv 1(mod \ 2)$, we obtain that $$CF_2(Dic_{4n})=CF_2(D_{2n})=2n.$$ 
The same reasoning is used if $n$ is an odd number. The only change is that, in this case, we can find other cyclic factorizations of $Dic_{4n}$ besides those corresponding to the ones of $D_{2n}$. Let $K$ be one of the subgroups contained in the set $\lbrace \langle a^{i-1}\gamma\rangle \ | \ i=\overline{1,n}\rbrace$. Again, if $H$ is a nontrivial cyclic subgroup of odd order contained in $\mathbb{Z}_{2n}$, then $(H,K)$ is a cyclic factorization of $Dic_{4n}$ if $|H|=n$. According to the subgroup structure of the dicyclic group, there is only one such subgroup isomorphic to $\mathbb{Z}_n$. It follows that
$$CF_2(Dic_{4n})=CF_2(D_{2n})+2n=4n,$$ and our proof is complete.         
\hfill\rule{1,5mm}{1,5mm}\\

The last class of groups we study is formed by the generalized dicyclic groups $Dic_{4n}(A)=\langle A, \gamma \ | \ \gamma^4=e, \gamma^2\in A\setminus \lbrace e\rbrace, g^{\gamma}=g^{-1}, \forall g\in A\rangle$, where $A\cong \mathbb{Z}_2\times \mathbb{Z}_n$. Let $a$ and $b$ be the generators of the cyclic groups $\mathbb{Z}_n$ and $\mathbb{Z}_2$, respectively. Since $\gamma^2\in A\setminus \lbrace e\rbrace$ and $\gamma^4=e$, we infer that $\gamma^2\in \lbrace a^{\frac{n}{2}}, b, a^\frac{n}{2}b \rbrace$. The subgroup strucute and explicit formulas for the (cyclic) subgroup commutativity degree of $Dic_{4n}(A)$ are provided by \cite{15} and our final purpose in this paper is to compute its cyclic factorization number.\\

\noindent\textbf{Theorem 3.9.} \textit{Let $A=\mathbb{Z}_2\times \mathbb{Z}_n$ be an abelian group, where $n=2^mm'$, $m\in\mathbb{N}^*$ and $m'$ is a positive odd integer. Then, \\ \\
(i) \ if $m=1, m'\not=1$ and $\gamma^2 \in \lbrace a^{\frac{n}{2}}, b, a^\frac{n}{2}b \rbrace$ or if $m\geqslant 2, m'$ is any odd number and $\gamma^2\in\lbrace b,a^{\frac{n}{2}}\rbrace$, the cyclic factorization number of the generalized dicyclic group $Dic_{4n}(A)$ is 
$$CF_2(Dic_{4n}(A))=4n.$$\\
(ii) \ if $m\geqslant 2, m'\not=1$ and $\gamma^2=a^{\frac{n}{2}}$, the cyclic factorization number of the generalized dicylic group $Dic_{4n}(A)$ is
$$CF_2(Dic_{4n}(A))=0.$$\\
(iii) \ if $m\geqslant 2, m'=1$ and $\gamma^2=a^{\frac{n}{2}}$, the generalized dicyclic group $Dic_{4n}(A)$ is isomorphic to the direct product $\mathbb{Z}_2\times Q_{2^{m+1}}$ and its cyclic factorization number is
$$CF_2(\mathbb{Z}_2\times Q_{2^{m+1}})=0.$$\\
(iv) if $m=m'=1$ and $\gamma^2\in \lbrace a^\frac{n}{2}, b, a^\frac{n}{2}b\rbrace$, the generalized dicyclic group $Dic_{4n}(A)$ is isomorphic to the abelian group $\mathbb{Z}_2^3$ and its cyclic factorization number is $$CF_2(\mathbb{Z}_2^3)=0.$$}
\\

\noindent\textbf{Proof.} The cyclic subgroups of $Dic_{4n}(A)$ are those contained in the abelian group $A$ and the subgroups contained in the set $C_1=\lbrace\langle a^{i-1}\gamma\rangle \ | \ i=\overline{1,n}\rbrace$ if $\gamma^2\in\lbrace b,a^{\frac{n}{2}}b\rbrace$, respectively in the set $C_2=\lbrace\langle (ab)^{i-1}\gamma\rangle \ | \ i=\overline{1,n}\rbrace$ if $\gamma^2=a^{\frac{n}{2}}$.

(i) Let $m\geqslant 2, m'$ a positive odd number and $\gamma^2\in\lbrace b, a^{\frac{n}{2}}b\rbrace$.  Since $|Dic_{4n}(A)|=4n$ and, in this case, $4n\geqslant 16$ we can not form a cyclic factorization using only cyclic subgroups contained in $C_1$ which are of order 4 and do not intersect trivially. Hence a cyclic factorization $(H,K)$ of $Dic_{4n}(A)$ is formed by a cyclic subgroup $H$ of $A$ and a cyclic subgroup $K$ belonging to $C_1$. This leads us to $|H|\geq n$ and since $A\cong \mathbb{Z}_2\times\mathbb{Z}_n$, it is clear that $|H|=n$. In other words, $H$ is a maximal cyclic subgroup of $A$. It is easy to see that $A$ has 2 maximal cyclic subgroups, namely $\langle a\rangle$ and $\langle ab\rangle$. We remark that the intersection of each one of these subgroups with any subgroup contained in $C_1$ is trivial. Hence, $CF_2(Dic_{4n}(A))=4n.$ If $m=1, m'\not=1$ and $\gamma^2 \in \lbrace b, a^\frac{n}{2}b \rbrace$, the reasoning is similar, the only change being that $A$ contains 3 maximal cyclic subgroups, namely $\langle a\rangle, \langle ab\rangle$ and $\langle a^2b\rangle$. It is easy to see that only two of these subgroups trivially intersect the $n$ subgroups belonging to $C_1$ as we choose $\gamma^2$ to be $b$ or $a^{\frac{n}{2}}b$, so we obtain again that $CF_2(Dic_{4n}(A))=4n.$ If $m=1, m'\not=1$ and $\gamma^2=a^{\frac{n}{2}}$, the cyclic factorizations of $Dic_{4n}(A)$ are formed by the cyclic subgroups contained in $C_2$ and two of the three maximal cyclic subgroups contained in $A$, namely $\langle ab\rangle$ and $\langle a^2b\rangle$. Therefore, in all the mentioned cases, the cyclic factorization number of the generalized dicyclic group $Dic_{4n}(A)$ is $$CF_2(Dic_{4n}(A))=4n.$$

\noindent(ii) Let $m\geqslant 2, m'\not=1$ and $\gamma^2=a^{\frac{n}{2}}$. The only possible cyclic factorizations $(H,K)$ of $Dic_{4n}(A)$ are formed by the factors $H\in L_1(A)$ such that $|H|=n$ and $K\in C_2$. We saw that $A$ has two maximal cyclic subgroups, but since their intersection with the $n$ cyclic subgroups contained in $C_2$ is $\lbrace e,\gamma^2\rbrace$, it is clear that $$CF_2(Dic_{4n}(A))=0.$$

\noindent(iii) Let $m\geqslant 2, m'=1$ and $\gamma^2=a^{\frac{n}{2}}$. For more details about the mentioned isomorphism, we refer the reader to \cite{15}. Denote by $x, y, b_1$ and $e_1$ the elements $a\langle \gamma^2\rangle, \gamma\langle \gamma^2\rangle, b\langle \gamma^2\rangle$ and $e\langle \gamma^2\rangle$, respectively. We get the following isomorphism 
$$\frac{\mathbb{Z}_2\times Q_{2^{m+1}}}{\langle \gamma^2\rangle}=\langle b_1\rangle\times\langle x, y \ | \ x^{2^m}=y^2=e_1, x^y=x^{-1}\rangle\cong \mathbb{Z}_2\times D_{2n'},$$ where $n'=\frac{n}{2}=2^{m-1}.$ This isomorphism implies that the group $\mathbb{Z}_2\times D_{2n'}$ contains the subgroup lattice of $A'=\langle b_1,x\rangle\cong \mathbb{Z}_2\times \mathbb{Z}_{n'}$ obtained by quotiening $A$ by $\langle \gamma^2 \rangle$. In the same paper, we saw that besides the cyclic subgroups contained in $A'$, the group $\mathbb{Z}_2\times D_{2n'}$ has other $2n'$ cyclic subgroups of order 2, namely $\langle x^{i-1}y\rangle$ and $\langle x^{i-1}yb_1\rangle$, where $i=\overline{1,n'}$. Since a factor $H$ of a possible cyclic factorization $(H,K)$ of $\mathbb{Z}_2\times D_{2n'}$ is one of the $2n'$ above mentioned cyclic subgroups of order 2, the other factor $K$ must be a cyclic subgroup contained in $A'$ such that $|K|\geq 2n'$. However, $A'\cong \mathbb{Z}_2\times \mathbb{Z}_{n'}$, so there is no such $K$. It follows that 	
$$CF_2(\mathbb{Z}_2\times D_{2n'})=CF_2(\mathbb{Z}_2\times D_{2^m})=0.$$
A nontrivial subgroup of $\mathbb{Z}_2\times Q_{2^{m+1}}$ contains $\langle \gamma^2\rangle=\langle a^{n'}\rangle$ or it is one of the minimal subgroups $\langle b\rangle$ or $\langle a^{n'}b\rangle$. Assume that $(H,K)$ is a cyclic factorization of $\mathbb{Z}_2\times Q_{2^{m+1}}$ such that $\langle \gamma^2\rangle\subseteq H$ and $\langle \gamma^2\rangle\subseteq K$. Then there exists a corresponding cyclic factorization of $\mathbb{Z}_2\times D_{2n'}$, which is false since we proved that $CF_2(\mathbb{Z}_2\times D_{2n'})=0$. It is clear that other possible cyclic factorization $(H,K)$ of $\mathbb{Z}_2\times Q_{2^{m+1}}$ is formed by a factor $H$ being $\langle b\rangle$ or $\langle a^{n'}b\rangle$ and a factor $K$ that contains $\langle \gamma^2\rangle$. Since $|H|=2$, it follows that the order of the cyclic subgroup $K$ is $2^{m+1}$ or $2^{m+2}$. However, the group $\mathbb{Z}_2\times Q_{2^{m+1}}$ do not have any element that could generate such a cyclic subgroup. Therefore, we proved that 
$$CF_2(\mathbb{Z}_2\times Q_{2^{m+1}})=0.$$        

\noindent(iv) Let $m=m'=1$ and $\gamma^2\in \lbrace a^{\frac{n}{2}}, b, a^{\frac{n}{2}}b\rbrace$. It is straightforward to check the isomorphism $Dic_{4n}(A)\cong\mathbb{Z}_2^3$. Also, by \textbf{Theorem 3.1.}, we have $$CF_2(\mathbb{Z}_2^3)=0,$$ as desired.     
\hfill\rule{1,5mm}{1,5mm}\\
  
\section{Conclusions and further research}

All our previous results show that the cyclic factorization number
constitutes an interesting computational aspect of finite groups.
It is clear that the study started here can successfully be
extended to other classes of finite groups. This will surely be
the subject of some further research.
\bigskip

We end our paper by formulating several open problems concerning
this topic.

\bigskip\noindent{\bf Problem 4.1.} Compute explicitly the cyclic
factorization number of a finite Zassenhaus metacyclic group (see
e.g. \cite{4}, I).

\bigskip\noindent{\bf Problem 4.2.} Determine the class of finite groups/$p$-groups
$G$ for which $CF_2(G)\neq 0$.

\bigskip\noindent{\bf Problem 4.3.} Study other connections between
the cyclic factorization number and the factorization number of a
finite group.

\bigskip\noindent{\bf Problem 4.4.} Given a finite non-nilpotent group
$G$, describe the manner in which $CF_2(G)$ depends on the cyclic
factorization numbers of the Sylow subgroups of $G$.

\bigskip\noindent{\bf Problem 4.5.} What can be said about two finite
groups of order $n$ having the same cyclic factorization number\,?

\vspace*{3ex}
\small

\begin{minipage}[t]{7cm}
Marius T\u arn\u auceanu \\
Faculty of  Mathematics \\
"Al.I. Cuza" University \\
Ia\c si, Romania \\
e-mail: {\tt tarnauc@uaic.ro}
\end{minipage}
\hfill
\begin{minipage}[t]{7cm}
Mihai-Silviu Lazorec \\
Faculty of  Mathematics \\
"Al.I. Cuza" University \\
Ia\c si, Romania \\
e-mail: {\tt Mihai.Lazorec@student.uaic.ro}
\end{minipage}
\end{document}